# A Catalog of Facially Complete Graphs

James Tilley (Bedford Corners, NY; jimtilley@optonline.net)
Stan Wagon (Macalester College, St. Paul, MN; wagon@macalester.edu)
Eric Weisstein (Wolfram Research, Inc., Champaign, IL; eww@wolfram.com)

**Abstract.** Considering regions in a map to be adjacent when they have nonempty intersection (as opposed to the traditional view requiring intersection in a linear segment) leads to the concept of a *facially complete* graph: a plane graph that becomes complete when edges are added between every two vertices that lie on a face. Here we present a complete catalog of facially complete graphs: they fall into seven types. A consequence is that if $q$ is the size of the largest face in a plane graph $G$ that is facially complete, then $G$ has at most $\lfloor \frac{3}{2} q \rfloor$ vertices. This bound was known, but our proof is completely different from the 1998 approach of Chen, Grigni, and Papadimitriou. Our method also yields a count of the 2-connected facially complete graphs with $n$ vertices. We also show that if a plane graph has at most two faces of size 4 and no larger face, then the addition of both diagonals to each 4-face leads to a graph that is 5-colorable.

## 1. Introduction

The classic view of the map-coloring problem considers two regions (taken to be closed sets) as adjacent when their boundaries intersect in a segment of nonzero length, as opposed to merely a single point. But single-point intersections do arise in real-world maps: there is the famous Four Corners point (a *quadripoint*) in the USA, where Colorado, Utah, Arizona, and New Mexico meet, and many other examples [6, 19]. Indeed, in the Catania province of Sicily there are 10 towns that meet at Mt. Etna's summit (Figure 1). Because Bronte surrounds Maletto and so strikes Mount Etna twice, it can be considered as 11 regions with a common point.

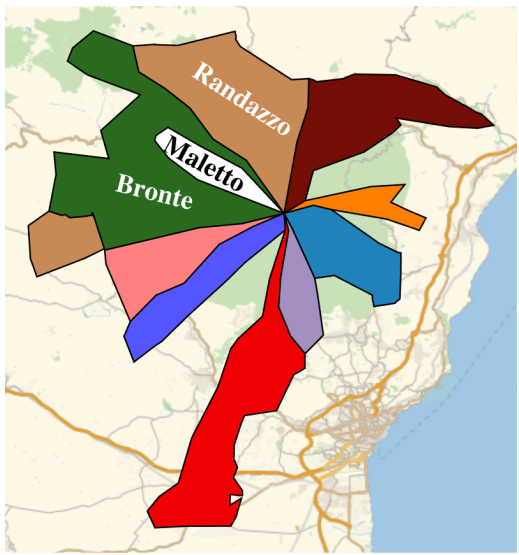

Figure 1. A map of part of the province of Catania in Sicily, showing an 11-fold meeting point of towns at the summit of Mount Etna. Only ten towns are involved because Bronte strikes the corner point twice, fully surrounding Maletto. The township of Randazzo is disconnected.

It is natural to require that regions with any nonempty intersection, even just a single point, be given different colors. It is not hard to construct a map with six regions that is complete (all pairs are adjacent) when 4-corner adjacency is allowed; see Figure 2, which is a simple modification of the four-corner layout to get a complete map with three quadri-points. The triangular representation is an alternate view. The adjacency graph at far right is $K_6$.

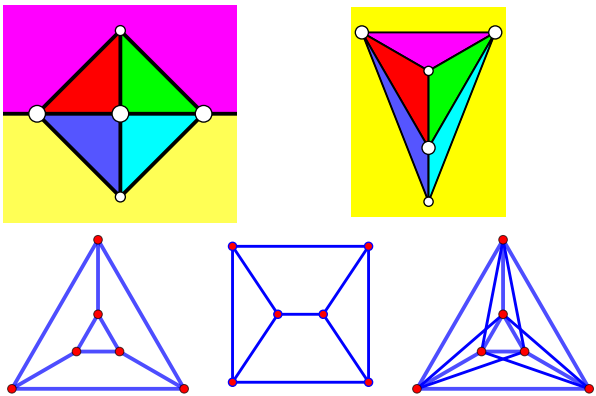

Figure 2. Top: Two views of a 5-node, 6-region map in which every two regions share a point; there are three 4-corner points. Bottom: The dual graph of the basic map, which is a 3-prism, shown in two drawings and followed by the graph (isomorphic to $K_6$) with the edges from the 4-corner adjacencies added.

For a general conjecture, we must deal with the map made up of sectors of a circle meeting at the center. Such a map shows that the chromatic number is unbounded, and therefore any coloring conjecture must involve a chromatic bound that is a function of $q$, the maximum number of countries that can share a point. This concept was first investigated in 1969 by Ore and Plummer [13], though Ringel [15] had considered a closely related problem a few years earlier.

As always, it is more convenient to work with graphs. A *cyclic coloring* of a plane graph is an assignment of colors to vertices such that any two vertices that share a face receive different colors. The face condition is the same as asking that regions in the dual map that share a common point get colored differently. The largest face in a plane graph is usually denoted by $\Delta^*(G)$, but here we will use $q$ throughout for the largest face size in $G$, where it will be clear from the context what $G$ is.

**Definition.** The *facial closure* $\overline{G}$ of a plane graph $G$ is defined by adding to $G$ any edge $u \leftrightarrow v$ not in $G$ whenever $u$ and $v$ lie on a face of $G$. A plane graph $G$ is *facially complete* (FC) if $\overline{G}$ is a complete graph.

Plummer and Toft [14] introduced $\chi_c$ for the cyclic chromatic number of a graph. We prefer to work with $\overline{G}$; $\chi_c(G)$ is the same as $\chi(\overline{G})$. Ore and Plummer [13] proved that $\chi(\overline{G}) \le 2\,q$ and presented the type 1 graphs—subdivided prisms—of §2: those graphs are facially complete and have $\lfloor \frac{3}{2} q \rfloor$ vertices, which means that the best possible conjecture would be $\chi(\overline{G}) \le \lfloor \frac{3}{2} q \rfloor$. The upper bound of $2\,q$ was improved by Sanders and Zhao [16] to $\chi(\overline{G}) \le \lceil \frac{5}{3} q \rceil$. The prism-based family and the cyclic coloring conjecture were also implicit in the work of Borodin [1]. So the current main question in this area is the following conjecture (see [8, §2.5]). We use $W(q)$ for $\lfloor \frac{3}{2} q \rfloor$.

**Cyclic Coloring Conjecture.** If $G$ is a plane graph, then $\chi(\overline{G}) \le W(q)$.

The $q = 3$ case of this conjecture is true as it is equivalent to the four-color theorem. The $q = 4$ case was proved by Borodin [1, 3] and the $q = 6$ case was proved by Hebdige and Krâl [5].

The reason $W(q)$ is the conjectured bound on $\chi(\overline{G})$ goes beyond the fact that there exists a facially complete graph of that size. It is known that that there is no larger facially complete graph.

**Theorem.** A plane graph that is facially complete has at most $W(q)$ vertices.

This result is stated at the end of the 1998 paper of Chen et al [4], but their investigation focuses on maps. In this paper we stay in the realm of graphs and present a straightforward proof of the theorem. Our approach is quite different from that of [4]. Further, we categorize the collection of all facially complete plane graphs using seven types. The first five types are the 2-connected cases and our description allows us to present a method that yields the number of 2-connected facially complete graphs with $n$ vertices.

Computations with graphs were used heavily in our investigations. We used *Mathematica*, applying various graph theory functions to the graphs in the `GraphData` database.

We are grateful to Joan Hutchinson, Paul Kainen, Brendan McKay, and John Watkins for helpful discussions about coloring and planar graph embedding.

# 2. The FC Catalog

An *outerplanar graph* is traditionally defined as a graph that has a drawing in which all vertices belong to the exterior face. Here we use the term for a plane graph having all vertices on one face; that is, we use it to describe an embedding as opposed to the abstract graph. The following gives a catalog of FC graphs. In §3 we will prove that every FC plane graph belongs to one of these families. We use the term *subdivider* for a vertex of degree 2 and we use $\kappa$ to denote the vertex connectivity of a graph.

**The FC Catalog.** Each plane graph in the following seven families is FC. For graphs in the first five families $\kappa \geq 2$; in the sixth family $\kappa = 1$, and in the last $\kappa = 0$.

*Type 1.* A subdivided 3-prism: subdividers are added to the edges connecting one triangle to another (Figure 3).

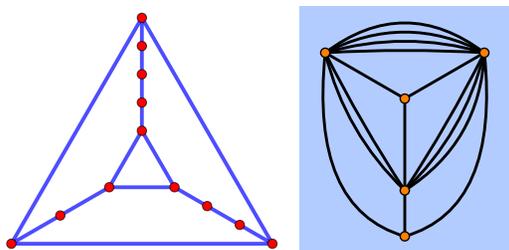

Figure 3. The 3-prism with 1, 2, and 3 subdividers; the corresponding 5-node map is at right with the exterior face in the graph represented by the bottom node of the map.

*Type 2.* A subdivided tetrahedral graph $K_4$: subdividers are added to the edges connecting a point to an exterior triangle (Figure 4). To avoid duplication with type 4, each connecting edge is to have at least one subdivider (otherwise the graph is a wheel).

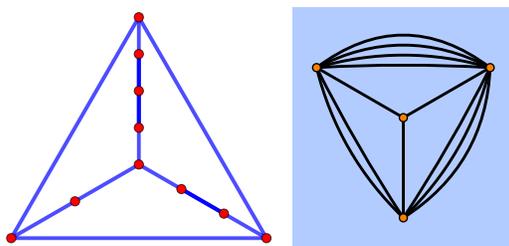

Figure 4. The tetrahedral graph $K_4$ with 1, 2, and 3 subdividers; the 4-node map is at right, with the exterior face of the graph represented by the central node of the map.

*Type 3.* A subdivided complete bipartite graph $K_{2,3}$: if the small part has vertices $\alpha_0$ and $\alpha_1$, then subdividers are added to the three disjoint paths from $\alpha_0$ to $\alpha_1$ (Figure 5). To avoid duplication with the wheels in type 4, each of the three paths must have at least two subdividers (otherwise the graph is a wheel with all but two radii deleted).

*Type 3.* A subdivided complete bipartite graph $K_{2,3}$. In order to avoid the graph being a wheel with all but two radii deleted (type 4), each subdivider in $K_{2,3}$ must have an adjacent subdivider (see Figure 5).

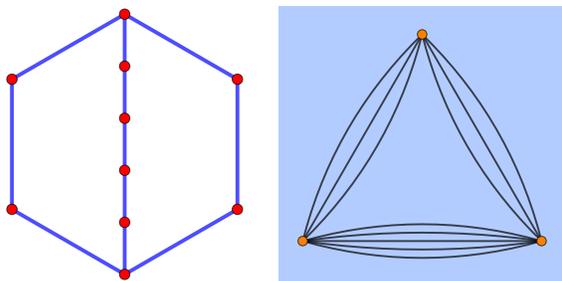

Figure 5. The subdivided $K_{2,3}$. The corresponding 3-node map is at right, with the top node representing the exterior face of the graph.

*Type 4.* The $n$-wheel where $n \geq 4$, with radial edges removed leaving a non-outerplanar graph. If the graph were outerplanar it would be a cycle with 0 or 1 diagonals (type 5) or a wheel with one radius (type 6).

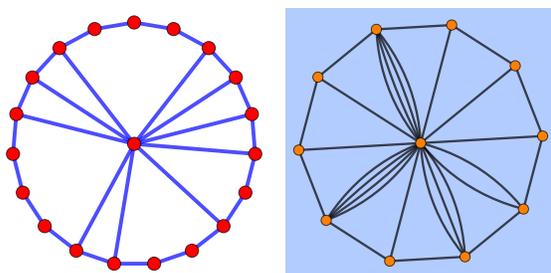

Figure 6. A 20-wheel with 9 radii removed. The dual map is at right, with the central node representing the graph's exterior face.

*Type 5.* A simple cycle with interior diagonals leading to a plane graph (Figure 7).

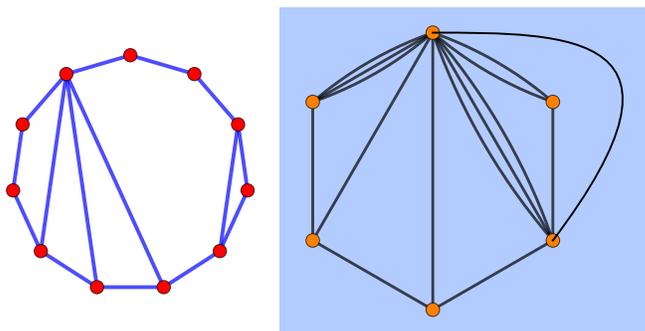

Figure 7. An 11-cycle with diagonals; the dual map is at right. The top node represents the graph's exterior face.

*Type 6.* A connected outerplanar graph with vertex-connectivity 1, or $K_1$. Alternatively, a connected cactus graph that is not a cycle. A *cactus graph* is one in which every block is an edge or a cycle (see Figure 8).

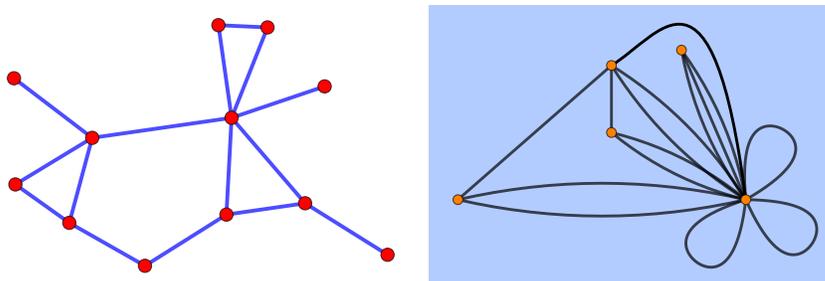

Figure 8. A 1-connected outerplanar graph, with its dual map at right (the lower right node is common to all 12 countries).

*Type 7.* A disconnected outerplanar graph; that is, disjoint copies of graphs of types 5 and 6.

The following two propositions present properties of the graphs in the catalog. We will prove in §3 that the catalog is complete, so these properties apply to all facially complete plane graphs. The extremal cases in (b) are of special interest; the type 1 family was presented by Borodin, and by Plummer and Ore as noted earlier. The type 2 family seems to have been overlooked. Note that the type 1 family in the maximal case for $q$ even is symmetric in that the number of subdividers is the same for each subdivided edge; this is not true when $q$ is odd. But for odd $q$, when the type 2 family is extremal, it is symmetric. The theorem in §1 follows from part (b) below once we prove that the catalog is complete.

**Proposition 1(a).** The seven types in the FC catalog are mutually exclusive.

**(b).** For all graphs in the FC catalog, $n \leq W(q)$. The only cases where $n = W(q)$ are as follows (see Figure 9):

• type 1, the subdivided 3-prism, with the subdividers distributed as $(a, a, a)$ or $(a, a, a + 1)$, where $a \geq 0$ (here $q \geq 4$ and $q$ is even in the first case and odd in the second);

• type 2, the subdivided tetrahedral graph, with subdivider distribution $(a, a, a)$, where $a \geq 1$ (this implies $q \geq 5$ and $q$ is odd);

• $K_4$, which is the 4-wheel in type 4 ($q = 3$).

**(c).** Except for $K_4$ and the $n$-wheels with $n$ even, which are 4-chromatic, the chromatic number of any graph in the catalog is at most 3.

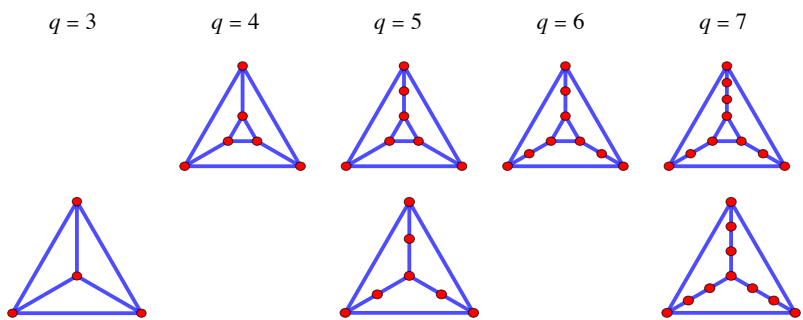

Figure 9. The graphs in the FC catalog (types 1 and 2, and $K_4$) having maximum vertex count $\lfloor \frac{3}{2} q \rfloor$.

*Proof (a)*. Types 6 and 7 are distinguished from each other and the others by connectivity ($K_1$ being a special case). A graph in type 5 is outerplanar while types 1–4 are not. Types 1, 2, and 3 have, respectively, 6, 4, and 2 vertices of degree 3. A wheel with four or more radii has a vertex of degree 4, and so cannot have type 1, 2, or 3. A wheel with two radii is potentially of type 3, but was excluded there, except that if the two radii involve neighbors on the rim, it is outerplanar. That case was excluded as it is of type 5. A wheel with three radii is potentially of type 2 when two of the radii involve rim neighbors, but was excluded from type 2.

*(b)*. For types 1–3, suppose the subdivider counts, in nondecreasing order, are $a$, $b$, $c$. In each case we can describe $q$ and $n$ in a way that yields the desired inequalities.

Type 1. $q = b + c + 4$ and $n = a + b + c + 6$. If $q$ is even, $W(q) - n = \frac{3}{2} q - n = \frac{1}{2} (b + c - 2 a)$. This last is nonnegative and is 0 only when $a = b = c$. If $q$ is odd, $W(q) - n = \frac{1}{2} (3 q - 1) - n = \frac{1}{2} (b + c - 2 a - 1)$. Because $b + c$ is odd, and $b + c \geq 2 a$, we have $b + c - 2 a \geq 1$. This proves the inequality. For equality, we must have $b + c = 2 a + 1$, from which the claimed form follows.

Type 2. $q = b + c + 3$ and $n = a + b + c + 4$. The proof is similar to the type 1 case. Details: If $q$ is even, $W(q) - n = \frac{3}{2} q - n = \frac{1}{2} (b + c - 2 a + 1)$. This last is positive. If $q$ is odd, $W(q) - n = \frac{1}{2} (3 q - 1) - n = \frac{1}{2} (b + c - 2 a)$. This is nonnegative and is 0 only if $a = b = c$.

Type 3. $q = b + c + 2$ and $n = a + b + c + 2 \leq 2 b + c + 2$. Here $2 (W(q) - n) \geq 2 \left( \frac{3 q - 1}{2} - n \right) = b + c - 2 a + 1 \geq 1$, so $W(q) - n \geq \frac{1}{2}$ and therefore $n < W(q)$.

Type 4. $n = q + 1 \leq W(q) - 1$, except when $G = K_4$. In that case $n = 4 = W(3)$.

Types 5, 6, and 7. $n = q \leq W(q) - 1$.

*(c)*. The only difficult case is the cycle with diagonals, but the 3-colorability of the vertices in the graph of a polygon triangulation is a well-known result in the context of the art gallery theorem (see [17]). The proof is by otectomy: any triangulation has an ear; cut it off; use induction. □

If 11 countries are allowed to meet at a point, as essentially happens for the towns in Sicily in Figure 1, a map such as the one in Figure 10 shows that $K_{16}$ can arise and so the map requires 16 colors. The map in Figure 10 is of the subdivided tetrahedral type; because 11 is odd there is a second, less symmetric, example arising from the subdivided prism with 16 vertices and subdivider pattern 4, 3, 3.

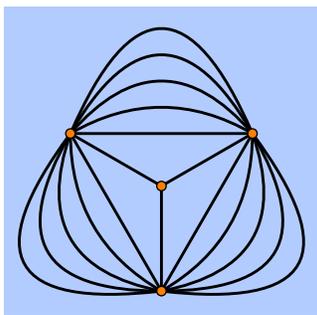

Figure 10. A Mt. Etna–inspired map with 16 countries having 11 of them meet at each of three points in such a way that every two regions share a point. The adjacency graph is therefore $K_{16}$ and 16 colors are needed.

We can determine the exact number of 2-connected FC graphs. There is an explicit formula in all cases except type 5, in which case there is an algebraic algorithm.

**Proposition 2.** The counts of the $n$-vertex graphs in the first five types in the FC catalog are as follows.

Type 1. $\left\lfloor \frac{(n-3)^2+6}{12} \right\rfloor$.    Type 2. $\left\lfloor \frac{(n-4)^2+6}{12} \right\rfloor$.    Type 3. $\left\lfloor \frac{(n-5)^2+6}{12} \right\rfloor$.

Type 4. $\frac{1}{4}\,(\mathrm{mod}(n-1,2)+3)\,2^{\lfloor (n-1)/2 \rfloor} + \frac{1}{2\,(n-1)}\left(\sum_{d \mid n-1} \phi(d)\,2^{(n-1)/d}\right) - 3$, where $\phi$ is Euler's $\phi$ function.

Type 5. For $n \geq 3$, the counts are $1, 2, 3, 9, 20, 75, 262, 1117, 4783, 21\,971, 102\,249, 489\,077, \ldots$ .

*Proof.* The type 1 count is the number of partitions of $n-6$ into three nonnegative parts [9]). The type 2 count is the number of partitions of $n-7$ into three nonnegative parts. It is $n-7$ because we must place at least one subdivider on each edge from the center, so the initial graph has 7 vertices. The type 3 formula counts the number of partitions of $n-8$ into three nonnegative parts. The use of $n-8$ is correct because the initial graph has 8 vertices. The type 4 count arises from a bracelet counting theorem [10], subtracting 3 for the outerplanar cases. The values for type 5 were generated by the method at [11]. □

These counts, and consequently the types in the catalog, were confirmed for $n \leq 8$ by using an algorithm for all distinct plane embeddings applied to *Mathematica*'s `GraphData` database. This resulted in counts of all 2-connected $n$-vertex plane that are facially complete that agree with the preceding proposition. The counts of 2-connected FC plane graphs, for $n = 3, 4, \ldots$ are $1, 3, 6, 15, 32, 94, 295, 1169, 4870, 22\,110, 102\,490, 489\,479, \ldots$ [12]. The counts excluding the outerplanar cases — that is, types 1, 2, 3, and 4 — are $0, 1, 3, 6, 12, 19, 33, 52, 87, 139, 241, 402, \ldots$ .

# 3. The Catalog is Complete

To prove that any FC plane graph has one of the types in the catalog, we need a lemma. This lemma applies to all graphs, not just FC ones. If $F$ is a face of a plane graph then $\hat{F}$ is the ordered set of vertices defining $F$ (possibly with duplications).

**Lemma.** Let $G$ be a 2-connected plane graph drawn with $q$-face $A$ as the exterior face, where $q \geq 4$ and $n \geq q + 2$. There exists a path $(\alpha_0, \rho_0, \rho_1, \ldots, \rho_k, \alpha_1)$ in $G$ where $\alpha_i \in \hat{A}$, $\alpha_0 \neq \alpha_1$, $\alpha_0$ and $\alpha_1$ are not consecutive vertices in the cycle bounding $A$, and each $\rho_i \notin \hat{A}$.

*Proof.* Let $\{a_i : 1 \leq i \leq q\}$ enumerate $\hat{A}$ in cyclic order around the boundary of $A$; always view indices as being reduced modulo $q$. We can partition the vertices not in $\hat{A}$ into connected components $\{C_i : 1 \leq i \leq r\}$ when the vertices in $\hat{A}$ are deleted; because $n \geq q + 2$, there is at least one component.

Because we seek a path, we can consider $G^-$, the graph obtained by contracting each $C_i$ to a vertex $\gamma_i$; the vertices $\gamma_i$ are independent in $G^-$. If any $\gamma_i$ has degree 1, 2-connectedness of $G$ would be violated. If some $\gamma_i$ has degree 3 or more, the fact that $q \geq 4$ means that the desired path would exist. If some $\gamma_i$ is adjacent to nonconsecutive vertices in $\hat{A}$, the desired path would exist. To complete the proof we will obtain a contradiction from the adjacency of each $\gamma_j$ to a pair $a_{k_j}, a_{k_j+1}$ (see Figure 11).

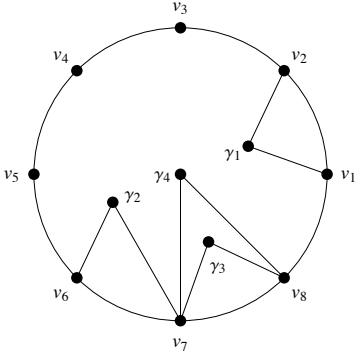

Figure 11. The contracted graph $G^-$ has a large interior face.

Each triple $\gamma_j, a_{k_j}, a_{k_j+1}$ forms either a triangular face or, with another $\gamma_i$, a 4-sided face. There must be at least one triangle because any quadrilateral has a $\gamma_m$ in its interior and $\gamma_m$, having degree 2, lies on one other face; repeat the argument until the second face is a triangle. Let $B$ be the face that is the complement of the exterior face and all the triangles or quadrilaterals just described. Delete all vertices and incident edges for any $\gamma_i$ that is not on $B$ and let $m$ be the number of $\gamma_i$ that remain; $m \geq 1$. Each remaining $\gamma_i$ contributes two edges to $B$ and it follows that

$|\hat{B}| = q - m + 2m = q + m > q$, contradiction. □

With the lemma in hand, we can prove that any FC graph is in the catalog. By Proposition 1(b) in §2 this implies that any facially complete graph has at most $W(q)$ vertices.

**Theorem.** Let $G$ be an FC plane graph. Then $G$ is one of the graphs in the FC catalog.

*Proof.* Assume that $G$ is drawn so that $|\hat{A}| = q$, where $A$ is the exterior face. We can assume $q \geq 4$, for if $q = 3$ then $\bar{G} = G$ and so $G$ is a complete graph; $K_1$ and $K_2$ are of type 6, $K_3$ has type 5, and $K_4$ has type 4. Assume first that $G$ is 2-connected.

$n = q$:  $G$ is a cycle with diagonals; type 5.

$n = q + 1$:  If two nonconsecutive vertices in $\hat{A}$ lie on an edge then the edge divides the complement of $A$ into two regions; neither region has a vertex in its interior, as such a vertex would lead to an independent pair in $\bar{G}$. Therefore there is no such edge and the additional vertex serves as the center of a wheel with radii deleted but at least two retained. This means $G$ is of type 4, or, if there are two radii involving neighbors in $\hat{A}$, $G$ is a cycle with one diagonal, which is type 5.

$n \geq q + 2$: Use the lemma and $q \geq 4$ to get a path $(\alpha_0, \eta_0, \ldots, \eta_1, \alpha_1)$ and let $P$ be the vertices of the path excluding its two ends; let $z = |P|$. Let $T$ be the vertices of one of the paths in $\hat{A}$ from $\alpha_0$ to $\alpha_1$, and let $S$ be be the same for the other path; exclude $\alpha_0$ and $\alpha_1$ from $S$ and $T$. Because $\alpha_0$ and $\alpha_1$ are not neighbors neither $T$ nor $S$ is empty. Let $t = |T|$ and $s = |S|$. Let $\beta_0 \in T$ and $\gamma_0 \in S$ be the neighbors of $\alpha_0$; let $\beta_1 \in T$ be the neighbor of $\alpha_1$. The path from the lemma divides the complement of $A$ into two open topological disks $B$ and $D$. There are no vertices in $B \cup D$, as any such would be separated from either $\beta_0$ or $\gamma_0$ in $\bar{G}$. Therefore, because $n \geq q + 2, z \geq 2$.

Let $T^+ = T \cup \{\alpha_0, \alpha_1\}$; let $P^+ = P \cup \{\alpha_0, \alpha_1\}$. Let $\epsilon$ count the edges that intersect the open set $B \cup D$. There are no edges in the $\epsilon$-count from $P^+$ to $P^+$ or from $\hat{A}$ to $\hat{A}$ as any such would separate a pair in $\bar{G}$. So any edge in $D$ must go from $T$ to $P$ and any edge in $B$ must go from $S$ to $P$.

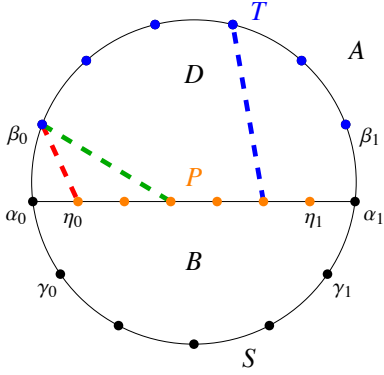

Figure 12. The orange points are $P$, the path given by the lemma. The blue points are $T$.

$\epsilon = 0$. $G$ is a subdivided $K_{2,3}$ and so has type 3 unless one of the bipartite edges has fewer than two subdividers, in which case it is a wheel, type 4.

$\epsilon = 1$. We may assume the edge is in $D$. Then $B$ is a face. If $t < 2$ then $|\hat{B}| > q$, contradiction, so $t \geq 2$. Any edge from $T \setminus \{\beta_0, \beta_1\}$ to $P$ (blue in Figure 12) separates a $\beta_i$ from the opposite $\eta_j$. And the same is true for an edge from $\beta_1$ to a vertex in $P$ other than $\eta_i$. So the edge in $D$ must be $\beta_0 \leftrightarrow \eta_0$ or $\beta_1 \leftrightarrow \eta_1$ (red in Figure 12). Therefore $G$ is in the tetrahedral case and has type 2.

$\epsilon = 2$. If both edges are in $D$ then, as in the $\epsilon = 1$ case, they must be $\beta_0 \leftrightarrow \eta_0$ and $\beta_1 \leftrightarrow \eta_1$ and $G$ is a subdivided prism, type 1. If there is one edge in each of $B$ and $D$ then again $G$ is type 1 because the $D$-edge can be taken to be $\beta_0 \leftrightarrow \eta_0$ and the $B$-edge is either $\gamma_0 \leftrightarrow \eta_0$ or $\gamma_1 \leftrightarrow \eta_1$, where $\gamma_1$ is the neighbor of $\alpha_1$ in $S$; $\gamma_1$ might equal $\gamma_0$. In the first case $\beta_0$ would be separated from $\eta_1$ in $\overline{G}$ and the second case is a graph of type 1.

$\epsilon \geq 3$. This cannot happen. We may assume there are at least two edges in $D$. If no edge is in $B$, then the argument of the $\epsilon = 1$ case shows that there are most two edges in $D$ and $\epsilon \leq 2$. Therefore there is an edge $\gamma_2 \leftrightarrow \eta_2$ in $B$ where $\eta_2 \in P$. We may assume that $\eta_2 \neq \eta_0$. Then any $D$-edge not containing $\eta_0$ separates $\eta_0$ from $\alpha_1$ in $\overline{G}$ and any $D$-edge containing $\eta_0$, but not equal to $\beta_0 \leftrightarrow \eta_0$ separates $\beta_0$ from $\eta_2$ in $\overline{G}$. So $D$ can have only the single edge $\beta_0 \leftrightarrow \eta_0$, a contradiction.

If $\kappa(G) = 1$, let $v$ be a cut-vertex of $G$. If necessary, modify the drawing of $G$ so that $v$ is on the exterior face. Then $v$ appears at least twice on the cycle defining the exterior face. The exterior face's cycle can be decomposed into at least two segments that correspond to simple cycles (in the case of a path, the segments will be singletons). Any vertex inside any one of these cycles would be separated in $\overline{G}$ from a vertex not on the cycle. Therefore the exterior face contains all vertices, and $G$ is outerplanar.

If $G$ is disconnected, then each connected component is FC and because there are at least two components, each component must be outerplanar and connected, and so must have type 5 or 6. □

# 4. Using Fewer Colors

In some cases we can show that $\chi(\overline{G})$ is strictly less than $W(q)$. The $q = 3$ case is nothing new as $G = \overline{G}$ and the cyclic coloring conjecture is equivalent to the four-color theorem. So suppose $q = 4$ and use the term *quad* for a 4-sided face. The 3-prism has three quads and its facial completion is 6-chromatic. We can show that any plane graph having at most two quads has a 5-colorable facial completion. The 5 is sharp because the 5-wheel has one quad and its facial completion is $K_5$, so five colors are necessary.

*Theorem.* If $G$ is a plane graph with $q \leq 4$ and at most two quads, then $\chi(\overline{G}) \leq 5$.

*Proof.* If there are no quads then $\overline{G} = G$ and is therefore 4-colorable. Suppose there is one quad. If both diagonals are edges in $G$, then again $G = \overline{G}$ and $\chi(\overline{G}) \leq 4$. If neither diagonal is an edge, add one. Let $G^+$ denote the plane graph that is $G$ with one diagonal in the quad. Use the four-color theorem to color $G^+$ with $1, 2, 3, 4$ and assume a triangle of the quad uses colors $1, 2, 3$ (see Figure 13). If color 4 is used at the other quad vertex, we are done. If color 3 is used at the fourth vertex, replace it with 5 to get a valid 5-coloring of $\overline{G}$.

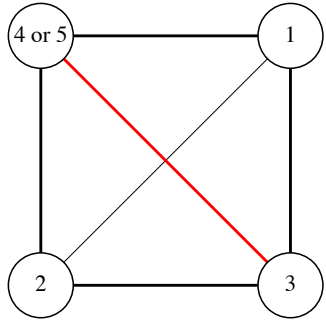

Figure 13. When there is one face of size 4, five colors suffice.

Suppose there are two quads, $A$ and $B$. If they have a common vertex, color each using the preceding case, making sure that the diagonals of the quads are such that if color 5 is used, it is used at a common vertex. If $A$ and $B$ are disjoint, add a diagonal to each and 4-color. This coloring respects any edges between the quads. Now use color 5 in each quad as in the one-quad case. But if the two vertices that get color 5 are on an edge, more work must be done to avoid using color 6. Note that at least one of the 16 possible edges from $A$ to $B$ is missing; for if they were all included, there would be a planar $K_{4,4}$ and therefore a planar $K_{3,3}$, which is impossible. Suppose one of these missing edges is $\alpha \leftrightarrow \beta$. Then color $A$ and $B$ so that, if color 5 is used, it is used only on $\alpha$ or $\beta$ or both. That is, use the one-quad method with diagonals chosen to be disjoint from $\{\alpha, \beta\}$. This is a legal 5-coloring of $\overline{G}$. □

# 5. Conclusion

The cyclic coloring conjecture of Plummer and Toft is an important extension of the four-color theorem and, like that theorem, it arises as a real-world map-coloring problem. The conjectured bound $W(q)$ bound on $\chi(\overline{G})$ is used because it appeared that there was no larger facially complete graphs. Our relatively simple proof that $W(q)$ bounds the size of an FC graph clarifies the result and our approach yields the fine detail about the family of all facially complete graphs.

Here are two natural questions about facial completions. There are several algorithms that can 4-color a planar graph. The proof yields one, though that would not be easy to implement. One can implement a random algorithm based on Kempe chains that apparently works with very high probability [7]. And there is also a deterministic algorithm that works well in practice [18]. Is there a reasonable algorithm that will succeed in coloring $\overline{G}$ in $W(q)$ colors?

Another question arises in the $q = 4$ case by trying to understand when the facial completion yields a nonplanar graph. It is not hard to see that the facial completion of the adjacency graph of the map of the US states is not planar. But the facial completion of a square is $K_4$, a planar graph. Is there a characterization or an algorithm that will tell when $\overline{G}$ is planar?

## References

1. O. V. Borodin, 1984. Solution of the Ringel problem on vertex-face coloring of planar graphs and coloring of 1-planar graphs, *Metody Diskret. Analiz.* 41 12–26 (in Russian).

2. O. V. Borodin, 2012. Coloring of plane graphs, A survey, *Disc. Math.*, 313:4, 517–539.

3. O. V. Borodin, 1995. A new proof of the 6 color theorem, *J. Graph Th.* 19, 507–521.

4. Z.-Z. Chen, M. Grigni, and C. H. Papadimitriou, 1998. Planar map graphs, *Proc. 30th ACM Symp. on Th. of Computing*, 514–523.

5. M. Hebdige and D. Krâl, 2016. Third case of the cyclic coloring conjecture, *SIAM J. Disc. Math.*, 30 525–548.

6. G. Hurst, 2016. Beyond Four Corners, USA, Wolfram Community, ⟨https://community.wolfram.com/groups/-/m/t/932548⟩.

7. J, P. Hutchinson and S. Wagon, 1998. Kempe revisited, *Amer. Math. Monthly,* 105, 170–174.

8. T. R. Jensen and B. Toft, 1995. *Graph Coloring Problems*, New York, Wiley.

9. Online Encyclopedia of Infinite Sequences, Partitions of *n* into at most 3 parts, ⟨https://oeis.org/A001399⟩.

10. Online Encyclopedia of Infinite Sequences, Bracelets with *n* beads, ⟨https://oeis.org/A000029⟩.

11. Online Encyclopedia of Infinite Sequences, Dissections of a polygon by nonintersecting diagonals, ⟨https://oeis.org/A001004⟩.

12. Online Encyclopedia of Infinite Sequences, Numbers of facially complete 2-connected planar embeddings, ⟨https://oeis.org/A375617⟩.

13. O. Ore and M. D. Plummer, 1969. Cyclic coloration of plane graphs, *Recent Progress in Combinatorics; Proc. 3rd Waterloo Conference on Combinatorics*, 287–293.

14. M. D. Plummer and B. Toft, 1987. Cyclic coloration of 3-polytopes, *J. Graph Th.* 11, 507–515.

15. G. Ringel, 1965. Ein Sechsfarbenproblem auf der Kugel, *Abhandlungen aus dem Mathematischen Seminar der Universität Hamburg*, 29 (1–2), 107–117.

16. D. P. Sanders and Y. Zhao, 2001. A new bound on the cyclic chromatic number, *J. Comb. Th. Ser. B*, 83, 102–111.

17. J. O'Rourke, 1987. *Art Gallery Theorems and Algorithms*, Oxford Univ. Pr.

18. J. Tilley, 2017. *D*-resolvability of vertices in planar graphs, *J. Graph Algorithms Appl,* 21(4), 649–661.

19. Wikipedia, Quadripoint, ⟨https://en.wikipedia.org/wiki/Quadripoint⟩.